\documentclass{radioeng}

\usepackage[czech,english]{babel}
\usepackage[cp1250]{inputenc} 
\usepackage{graphics}

\usepackage{amssymb,amsmath}
\usepackage{bm}
\usepackage{german}
\usepackage{graphicx}
\usepackage{cite}

\newcommand{\ZZ}{\mathbb{Z}}

\newcommand{\RR}{\mathbb{R}}

\begin{document}


\setcounter{page}{11}

\rengHeader{}{}{}{K. BITTNER, H.G. BRACHTENDORF, ADAPTIVE MULTI-RATE WAVELET METHOD FOR CIRCUIT SIMULATION}

\rengTitle{Adaptive Multi-rate Wavelet Method for Circuit Simulation}

\rengNames{Kai BITTNER$^1$, Hans Georg BRACHTENDORF$^1$}

\rengAffil{$^1$University of Applied Sciences Upper Austria, Softwarepark 11, A-4232 Hagenberg, Austria}

\rengMail{Kai.Bittner@fh-hagenberg.at, Hans-Georg.Brachtendorf@fh-hagenberg.at}

\begin{multicols}{2}

\begin{rengAbstract}
In this paper a new adaptive algorithm for
multi-rate circuit simulation encountered in the design of RF circuits 
based on spline wavelets is presented.
The circuit ordinary differential equations are first rewritten by a system of
(multi-rate) partial differential equations (MPDEs) in order to decouple the different time scales. 
Second, a semi-discretization by Rothe's method of the MPDEs
results in a system of differential algebraic equations (DAEs) with periodic boundary conditions.
These boundary value problems are solved by a Galerkin discretization using 
spline functions. An adaptive spline grid is generated, using spline wavelets for non-uniform grids.
Moreover the instantaneous frequency is chosen adaptively to guarantee a smooth envelope resulting in
large time steps and therefore high run time efficiency.
Numerical tests on circuits exhibiting multi-rate behavior including mixers and PLL conclude the paper.
\end{rengAbstract}

\rengKeywords{RF circuits, circuit simulation, multirate simulation, envelope simulation, splines, wavelets.}


\rengSection{Introduction}

Widely separated time-scales occur in many
radio-frequency (RF) circuits such as mixers, oscillators, PLLs, etc., 
making the analysis with standard
numerical methods difficult and costly. Low frequency or baseband signals and high
frequency carrier signals often occur in the same circuit,
enforcing very small time-steps over a long time-period
in the computation of the numerical solution. The occurrence of widely separated
time-scales is also referred to as a multi-scale or multi-rate problem.
Hence, classical
numerical techniques, as transient analysis, result into prohibitively long run-times. 

One method to circumvent
the bottleneck is to reformulate the ordinary circuit DAEs in a system of partial DAEs (multi-rate PDAE).
The method was first presented in \cite{Bra94,BWL+96}, specializing to multi-rate circuits 
and systems in steady-state.
In \cite{NL96,Roy97} the technique has been firstly applied to 
initial transient simulation of
driven circuits with a-priori known frequencies. In 
\cite{Bra97,BL98a,Bra2001} a generalization
has been derived for circuits with a-priori unknown or 
time-varying frequencies. These techniques opened the door
to multi-rate techniques with frequency modulated signal sources
or autonomous circuits such as oscillators with a priori unknown fundamental
frequency \cite{BL98a,BL98,Bra2001}.

Meanwhile, numerous articles on this topic 
have been published, often
specializing the technique for applications, e.g.\  
chirp signals \cite{Roy01}. The following
list on multi-rate publications is therefore not comprehensive 
\cite{BWL97,Lan02,PG02,Roy01,MRCHD04,OP07,OP09,OP11,RMM+11}.

In recent time the expansion of the signal waveforms by a wavelet or spline basis
instead of trigonometric basis functions became increasingly popular
\cite{SN03,SGN07,DCB05,Dau05,Bra2009,BiDau10b, BiDau10a}.
Tri\-go\-no\-metric basis functions are badly suited for pulse-shaped signals due to a slow convergence
of the trigonometric series and the Gibbs' phenomenon, whereas spline basis
functions with their compact support enable an adaptive mesh with a local refinement and hence smaller system size.
It has been shown in \cite{Bra2009} that the spectra can be easily calculated from a cubic spline
basis which is very important for the acceptance of the approach from the community of circuit
designers.

Essential for the success of such methods is that the spline or wavelet expansion can be adapted to the particular
signal shape in an circuit simulation example. First steps have been done in \cite{BiDau10b, BiDau10a}, where an adaptive 
spline wavelet approach has been developed for  transient analysis, i.e., the solution of initial value problems.
Here we present an extension of this approach to multi-rate problems, where advantages of the approach can be fully explored.
Sect.~\ref{multi_rate} will give a short introduction to the multi-rate circuit simulation problem. In Sect.~\ref{semid}--\ref{pss}
we describe the discretization of the multi-rate circuit equations based on a spline expansion for the semi discretized
problem.  Sect.~\ref{adaptive} introduces the scheme for wavelet based adaptive grid generation. Numerical tests in 
Sect.~\ref{numtest} show the performance of the method.

\rengSection{The multi-rate circuit simulation problem\\ \label{multi_rate}}

We consider circuit equations in the charge/flux oriented
modified nodal analysis (MNA) formulation, which yields a mathematical model
in the form of a system of differential-algebraic equations
(DAEs):
\begin{equation}
  \label{eq_MNA_charge}
  \tfrac{d}{dt}q\big(x(t)\big)
        + g\big(x(t)\big) + s(t) = 0.
\end{equation}

\noindent Here $x(t)\in\RR^n$ is the vector of node potentials and specific branch
currents and $q(x)\in\RR^n$ is the vector of charges and fluxes. The vector
$g(x)\in\RR^n$ comprises static contributions, while $s(t)\in\RR^n$ contains the
contributions of independent sources. The DAEs in
(\ref{eq_MNA_charge}) are usually solved by integration formulas for
stiff systems.

To separate the different time scales the problem is
reformulated as a partial differential-algebraic equation (PDAE) as
\begin{equation}\label{multirate}
\tfrac{\partial}{\partial \tau} q\big(\hat{x}(\tau,t)\big)
+\omega(\tau)\,\tfrac{\partial}{\partial t} q\big(\hat{x}(\tau,t)\big)
+g\big(\hat{x}(\tau,t)\big)+\hat{s}\big(\tau,t\big)=0
\end{equation}

\noindent where $\omega(\tau)$ is an estimate of the (scaled) angular frequency.
The bivariate function $\hat{x}(\tau,t)$ is related to the univariate solution $x(t)$ of 
(\ref{eq_MNA_charge}) as follows.
For any solution $\hat{x}(\tau,t)$ of (\ref{multirate}) 
we get by
\begin{equation}\label{character}
x_\theta(t) = \hat{x}\big(t,\Omega_\theta(t)\big),\quad\Omega_\theta(t)=\theta+\int_0^t \omega(s)\,ds
\end{equation}

\noindent a solution of 
\begin{equation}
  \label{singlerate}
  \tfrac{d}{dt}q\big(x(t)\big)+ g\big(x(t)\big) = \hat{s}\big(t,\Omega_\theta(t)\big).
\end{equation}

\noindent Thus, if we choose $\hat{s}$ such that
\begin{equation}\label{char_source}
s(t) = \hat{s}\big(t,\Omega_0(t)\big)
\end{equation}

\noindent then the solution of (\ref{multirate})  provides also a solution of (\ref{eq_MNA_charge}).

Although the formulation (\ref{multirate}) is valid for any circuit,
it offers a more efficient solution only for certain types of
problems. This is the case if $\hat{x}(\tau,t)$ is periodic in $t$
and smooth with respect to $\tau$. Then, a semi-discretization with
respect to $\tau$ can be done resulting in a relative small number of
 periodic boundary problems
in $t$ for only a few discretization points $\tau_\ell$.
In the sequel we
will consider (\ref{multirate}) with periodicity conditions in
$t$, i.e., $\hat{x}(\tau,t)=\hat{x}(\tau,t+P)$ and suitable
initial conditions $\hat{x}(0,t)=X_0(t)$.
Here $P$ is an arbitrary but fixed period length, which results in a scaling of $\hat{x}(\tau,t)$ and $\omega(\tau)$
as we will see in the next section. Although
$P>0$ can be chosen arbitrarily only a few choices are of practical interest.
These are essentially $P=1$ or $P=2\pi$ if one wants to work in a fixed periodic setting
or one choses $P$ as the period of the carrier so that the scaling of $\hat{x}(\tau,t)$
corresponds more to the physical setting. 

As we have shown in \cite[Lemma~3.2]{BiBra12b} the additional (scaled) parameter 
$\omega(\tau)$ can be used to adapt to the instantaneous frequency $f(\tau)$ of the carrier signal. If the instantaneous frequency is known (e.g. through the sources)
one can set $\omega(\tau)=f(\tau)\,P$, to obtain optimal smoothness with respect to $\tau$. 
If the instantaneous frequency is not known in advanced we consider $\omega(\tau)$ as an additional unknown, which is
determined by the additional smoothness condition
\begin{equation}\label{min_cond}
\int_0^P \Big|\tfrac{\partial}{\partial\tau}\hat{x}(\tau,t)\Big|^2\,dt\ \to \min
\end{equation}

\noindent(see \cite{BiBra12b} for details).

\rengSection{Semi-discretization\label{semid}}

We discretize (\ref{multirate}) with respect to $\tau$ (Rothe's method) using Gear's BDF\footnote{Other multistep
method (e.g.~trapezoidal rule) can be used, too. However, for simplicity of the representation we consider
only the BDF-method, which was also used for the implementation.}
method of order $s$ and obtain
\begin{align}\nonumber 
\sum_{i=0}^s \alpha^k_i q\big(X_{k-i}(t)\big) + 
\omega_k\frac{d}{dt}q\big(X_k(t)\big)&+g\big(X_k(t)\big)+\hat{s}\big(\tau_k,t\big)=0\\
\label{rothe}&X_k(t) = X_k(t+P).
\end{align}
Here we have to determine an approximation $X_k(t)$ of
$\hat{x}(\tau_k,t)$ from known, approximative solutions
$X_{k-i}(t)$ at previous time steps $\tau_{k-i}$. $i=1,\ldots,s$.
As it is well known, the BDF coefficients $\alpha^k_i$ are chosen such that for any polynomial $p$ of degree up to $s$ 
the derivative is computed exactly, i.e.,
$$
\sum_{i=0}^s \alpha^k_i p(\tau_{k-i}) = p'(\tau_k). 
$$ 

Assuming that the $\tau_{k-i}$, $\omega_{k-i}$ and $X_{k-i}(t)$ are known and fixed
for $i>0$ we define
\begin{equation}\label{fkxt}
f_k(x,t) :=\alpha^k_0 q(x) +g(x)+\hat{s}\big(\tau_k,t\big)+\sum_{i=1}^s\alpha^k_i q\big(X_{k-i}(t)\big).
\end{equation}

\noindent Then $X_k$ is the solution of the periodic boundary value problem
\begin{align}
\label{envelope}
\omega\frac{d}{dt}q\big(x(t)) + f_k(x(t),t) = 0&\\
\label{pbc} x(t) = x(t+P).
\end{align}
The new problem (\ref{envelope}) is closely related to the original periodic steady state problem of the circuit, only modified
by the additional `source term' $\sum_{i=1}^s\alpha^k_i q\big(X_{k-i}(t)\big)$.

If $\omega=\omega_k\approx\omega(\tau_k)$ is not fixed in advance, then an additional condition is needed.
Following \cite{BiBra12b} we use the discretized version 
\begin{equation}
\label{min_cond2}
\int_0^P\big|X_k(t)-X_{k-1}(t)\big|^2\,dt\to\min
\end{equation}

\noindent of (\ref{min_cond}).

\rengSection{Multi-rate source term}

Except for autonomous systems we have to choose a multi-rate source term $\hat{s}(\tau,t)$, which satisfies
for the given source $s(t)$ the relation (\ref{char_source}), even if the frequency term $\omega(\tau)$ is not known
in advance. Therefore, we consider first a reference source term for constant $\omega(\tau)$.
We assume that there is a reasonable multi-rate
term $\tilde{s}(\tau,t)$ which satisfies the relation 
$$
s(t) = \tilde{s}(t,\tilde{\omega}\, t)
$$
 i.e., (\ref{char_source}) for a constant $\omega(\tau) = \tilde{\omega}$. For any non-constant $\omega(\tau)$
we can than choose 
$$
\hat{s}(\tau,t)= \tilde{s}(\tau, t+\sigma(\tau)) 
$$ 
where 
$$
\sigma(\tau)=\tilde{\omega}\;\tau-\Omega_0(\tau)=\int_0^\tau \tilde{\omega}-\omega(s)\, ds.
$$

If $\omega(\tau)$ is fixed, then $\sigma(\tau_k)$ needed for the computation of $f_k$ can be computed
exactly or at least with arbitrary accuracy. If  $\omega(\tau)$ is not known, then an approximation 
$\sigma_k\approx \sigma(\tau_k)$
has to computed together with $\omega_k$. We will use the quadrature formula
\begin{equation}\label{sigmak}
\sigma_k:=
\sum_{\ell=1}^k(\tau_\ell-\tau_{\ell-1})\Big((\tilde{\omega}-\omega_{\ell-1})(1-W)+(\tilde{\omega}-\omega_{\ell})W \Big).
\end{equation}

\noindent In particular, with weights 
$W=0$, $W=1$ and $W=\frac{1}{2}$ we obtain the left and right rectangular rule and the trapezoidal
rule, respectively.\footnote{We recommend the choice $W=\frac{1}{2}$ or $W=1$, since this leads to a more stable
behavior of the sequence $(\omega_k)$.} 
If $W>0$ we have to take into account $\sigma_k$ and thus $f_k(x,t)$ depends on the unknown 
$\omega=\omega_k$ and with (\ref{fkxt}) and (\ref{sigmak}) the derivative can be calculated as 
\begin{equation}\label{dfo}
\frac{d}{d\omega}f_k(x,t)=-(\tau_k-\tau_{k-1})\, W\;\frac{d}{dt} \tilde{s}(\tau_k, t+\sigma_k).
\end{equation} 

\rengSection{A spline Galerkin method for the periodic boundary value problems\\ \label{pss}}

In \cite{BiDau10a,BiDau10b}, an adaptive
spline wavelet method for
the initial value problem on the circuit equations
(\ref{eq_MNA_charge}) has been developed. For
the periodic problem we will modify this approach by using a
periodic basis. Then the periodic boundary conditions
(\ref{pbc}) are fulfilled automatically and has not to be enforced
by additional equations.

We want to approximate the solution of (\ref{envelope}) by a periodic 
spline function of order $m$. For given grid points $t_\ell$ with
$$
0 = t_0 < t_1 < \ldots < t_N = P
$$
we consider all $m-2$-times differentiable, $P$-periodic functions, which are
polynomials of degree less than $m$ at each subinterval
$(t_\ell,t_{\ell+1})$. The break points $t_\ell$ are called
spline knots, which are periodically extended by $t_{kN+\ell}=t_\ell+k\,P$.
Note, that the periodicity condition implies that the periodic spline
is a piecewise polynomial also with respect to the extended grid.

A stable and computational efficient basis for the linear space of spline
functions is constituted by the B-splines $N^m_\ell(t)$, which are uniquely 
determined by their minimal support $[t_\ell,t_{\ell+m}]$ and the partition of
unity $\sum_i N^m_i(t) = 1$ (for normalization). For more information
on spline functions and B-splines as well as for efficient computational 
methods we refer the reader to the detailed description in \cite{Boor78,Sch81}.   

To expand the periodic solution of (\ref{envelope}), we need periodic basis
functions, which we obtain by the periodized B-splines
$$
\varphi_\ell(t)=\sum_{k\in\ZZ} N^m_{\ell+kN}(t)=\sum_{k\in\ZZ} N^m_\ell(t-k\,P),\quad \ell=1,\ldots,N,
$$
which form a basis for the space of periodic spline functions.

Now we have to find a spline function 
\begin{equation}\label{xapprox}
\tilde{x}(t)=\sum_{i=1}^N c_k \varphi_k(t)
\end{equation}

\noindent which approximates the solution of (\ref{envelope}). Since $c_\ell\in\RR^n$ 
($n$ is the number of equations and unknowns in the circuit equations
(\ref{eq_MNA_charge}))
we have to determine $n\times N$ coefficients $c_{k,i}$.
Thus, we have to derive from (\ref{envelope}) $n\times N$
conditions, which ensure that $\tilde{x}$ indeed approximates the
solution of (\ref{envelope}). 
Taking the integral over $N$ subintervals we obtain
\begin{eqnarray}
\label{nonlinear}
\lefteqn{0=F_\ell(\bm{c},\omega):= \int_{\hat{t}_{\ell-1}}^{\hat{t}_\ell}
\omega\,\tfrac{d}{dt}q_k\big(\tilde{x}(t)\big) + f_k\big(\tilde{x}(t),t\big)\,dt}\\
\nonumber
&=& \omega\Big(q_k\big(\tilde{x}(\hat{t}_\ell)\big)-q_k\big(\tilde{x}(\hat{t}_{\ell-1})\big)\Big)
+ \int_{\hat{t}_{\ell-1}}^{\hat{t}_\ell}f_k\big(\tilde{x}(t),t\big)\,dt
\end{eqnarray}
for $\ell=1,\ldots,N$. This  $N$ vector valued equations determine
the vector coefficients $c_k$. Note, that the splitting points $\hat{t}_k$
do not coincide with the spline knots $t_k$, but they have to be chosen
in relation to the spline grid. In particular, we need due to periodicity 
that $\hat{t}_N-\hat{t}_0=P$.

The nonlinear system (\ref{nonlinear}) can be solved by Newton's method. If $\omega=\omega(\tau)$
is known in advance we have to solve the linear system $\bm{A\,d}_c=\bm{b}$ to determine the 
Newton correction $\bm{d}_c$ in
\begin{equation}\label{newtonupdate}
\bm{c}^{(k,j+1)}=\bm{c}^{(k,j)}-\bm{d}_c
\end{equation}

\noindent where $j$ is the Newton count.
Starting from a sufficiently good initial guess, e.g. $\bm{c}^{(k,0)}=\bm{c}^{(k-1)}$, one obtains usually
after sufficiently many steps a good approximation $\bm{c}^{(k)}=c^{(k,J)}$ of the solution of (\ref{nonlinear}).
In order to set up the linear system we have to compute the right hand side\footnote{For clarity of notation we 
use $\bm{c}=\bm{c}^{(k,j)}$ and $\omega=\omega _{k,j}$ in the sequel, if the relation is clear from the context.}
$\bm{b}=F(\bm{c},\omega)=\big(F_\ell(\bm{c},\omega)\big)_{\ell=1,\ldots,N}$ as well as the Jacobian 
(with respect to $\bm{c}$) $A=D_cF(\bm{c},\omega)$
with the matrix block entries
$$
\frac{\partial F_\ell(\bm{c},\omega)}{\partial c_k}
= \omega\Big(C\big(\tilde{x}(\hat{t}_\ell)\big)-C\big(\tilde{x}(\hat{t}_{\ell-1})\big)\Big)
+ \int\limits_{\hat{t}_{\ell-1}}^{\hat{t}_\ell} G\big(\tilde{x}(t),t\big)\,dt
$$
where the Jacobians $C=D_x q(x)$ and  $G=D_x g(x,t)$ are available in any circuit simulator, which is
based on the MNA formulation (\ref{eq_MNA_charge}). The integrals can be computed by a suitable quadrature
rule. In practical computations where the spline order is chosen usually as 3 or 4 one can use the Simpson rule
or the two point Gau{\ss} quadrature.

If $\omega_k = \omega(\tau_k)$ has to 
be determined during the computation, linearization of (\ref{nonlinear}) results in the underdetermined
system
\begin{equation}\label{newton_ud}
\bm{A\,d}_c+d_\omega\,\bm{z} =\bm{b}
\end{equation}

\noindent under condition (\ref{min_cond2}). Here $\bm{z} = \big( D_\omega F_\ell(\bm{c},\omega)\big)_{\ell=1,\ldots,N}$,
where with (\ref{dfo})
\begin{eqnarray*}
\lefteqn{D_\omega F_\ell(\bm{c},\omega)} \\
& =& q_k\big(\tilde{x}(\hat{t}_\ell)\big)-q_k\big(\tilde{x}(\hat{t}_{\ell-1})\big)  
            +\int_{\hat{t}_{\ell-1}}^{\hat{t}_\ell}\tfrac{d}{d\omega}f_k\big(\tilde{x}(t),t\big)\,dt\\
& =& q_k\big(\tilde{x}(\hat{t}_\ell)\big)-q_k\big(\tilde{x}(\hat{t}_{\ell-1})\big)\\ 
& &\quad      +~(\tau_k-\tau_{k-1}) 
  \Big(\tilde{s}(\tau_k, \hat{t}_{\ell-1}+\sigma_k)-\tilde{s}(\tau_k, \hat{t}_\ell+\sigma_k) \Big)\,W .
\end{eqnarray*}
For computational reasons we replace (\ref{min_cond2}) by a similar condition on the spline coefficients, namely
\begin{equation}\label{min_cond3}
\big\|\bm{c}^{(k)}-\bm{c}^{(k-1)}\big\|^2_2\to \min.
\end{equation} 

Following \cite{BiBra12b} we obtain the solution of (\ref{newton_ud}) which
satisfies (\ref{min_cond3}) by
$$
\bm{d}_c=\tilde{\bm{b}} - d_\omega\,\tilde{\bm{z}},
$$
where $\tilde{\bm{b}}=\bm{A}^{-1}\,\bm{b}$ and $\tilde{\bm{z}}=\bm{A}^{-1}\,\bm{z}$ are computed 
by solving the corresponding linear systems and 
$$
d_\omega= \frac{\tilde{\bm{z}}^T\big(\bm{c}^{(k,j-1)}-\bm{c}^{(k-1)}-\tilde{\bm{b}}\big)}{\tilde{\bm{z}}^T\tilde{\bm{z}}}
$$
where in addition to (\ref{newtonupdate}) $\omega_k$ is determined by the iteration
$$
\omega_{k,j+1}=\omega_{k,j}-d_\omega,\qquad \omega_{k,0}=\omega_{k-1}.
$$

\rengSection{Adaptive grid generation\label{adaptive}}

It is still open, how the spline knots $t_i$ are chosen. In particular, for functions with sharp
transients it is important that we put more grid points at the location of this refinement. Since 
those locations are often not known in advance, it is important that the grid is generated automatically
by an adaptive refinement. For this refinement we use spline wavelets for non uniform grids introduced
in \cite{Bit05b}. 

Starting on some coarse initial grid we solve (\ref{nonlinear}) by Newton's method and obtain a first approximation, 
which we denote as $X(t)$.
Now we apply a fast wavelet transform introduced in \cite{BiBra14a} and obtain
\begin{equation}\label{wave_exp}
X(t) = \tilde{X}(t) + \sum_k d_k \,\psi_k(t).
\end{equation}

\noindent The spline $\tilde{X}(t)$ is an approximation of $X(t)$ using only the even spline knots $t_{2k}$.
The $\psi_k$ are wavelets which carry detail information. In the setting of \cite{Bit05b} the wavelets
are compactly supported splines, with a chosen number of vanishing moments, which are localized
near $t_{2k+1}$. Thus, the coefficients $d_k$ are a measure of the local approximation error, which 
describes also the smoothness of $X(t)$. Therefore we insert additional knots in
the neighborhood of $t_{2k+1}$ if the coefficient $d_k$ exceeds a given threshold.
The spline representation of $X(t)$ for the new grid can then be computed efficiently by the Oslo algorithm
\cite{CLR80,LyMo86}. These refined spline expansion is used as a new initial guess for another Newton iteration.

This approach is repeated several times, with Newton tolerances and wavelet threshold reduced
for each refinement.  This leads to better and better approximations, which provide more and more information
for the grid refinements leading to an almost optimal grid. Furthermore, this approach is computational efficient,
since the first Newton iterations are computed on a relatively small grid. For the final refinements usually one 
Newton update is sufficient, since we have already an excellent initial guess. The whole process is stopped
if the Newton iteration after the last refinement gets below a given error bound. These error bound is indeed a very good estimate of the achieved error. For details on the refinement algorithm we refer to \cite{BiBra14a}.

In Sect.~\ref{pss} we suggested to use the solution $X_{k-1}$ of the previous time step $\tau_{k-1}$
as initial guess (predictor) for Newton's method. In fact any reasonable predictor should be based on solutions from
previous time steps such that the following approach should be used for other predictors, too. 
Due to the smoothness in $\tau$ the difference between
$X_{k-1}$ and $X_k$ should be small, and we can expect fast convergence of Newton's method.
However, the refinements in the adaptive method described above will increase the number of spline knots and 
thus the size of the problem.  Therefore, we will use an approximation of $X_{k-1}$ on a coarser grid as initial 
guess. An efficient approximation can be achieved by wavelet thresholding as follows. From the wavelet expansion 
(\ref{wave_exp}) of $X_{k-1}$ we remove terms with small wavelet coefficients $d_k$, i.e., an approximation is
obtained as
\begin{equation}\label{thresh}
Y(t) = \tilde{X}(t) + \sum_{|d_k|>\varepsilon} d_k \,\psi_k(t).
\end{equation}

\noindent
If the threshold $\varepsilon$ is chosen in the order of magnitude of the error tolerance of the method, than we can still
expect the initial guess $Y(t)$ to be as good as $X_\ell(t)$ itself, while the computational cost is reduced. Using the approach
from \cite{BiBra14a} removing $\psi_k(t)$ is equivalent to remove the knot $t_{2k+1}$. Thus, the method will
remove spline knots, which became dispensable due to a change in signal shape over several time steps $\tau_k$,
or which where not needed from the beginning due to the rough error estimate in the refinement described above.
This leads to an almost optimal grid for any given error tolerance. For more details on this grid coarsening
algorithm we refer to  \cite{BiBra14a}.
\pagebreak[4]
    
\rengSection{Numerical tests\label{numtest}}

The above algorithm is implemented in C++, using a C++ library for circuit description and device models.  

\rengSubsection{Gilbert mixer}

\begin{figure}
 \begin{center}
 \includegraphics[width=\columnwidth]{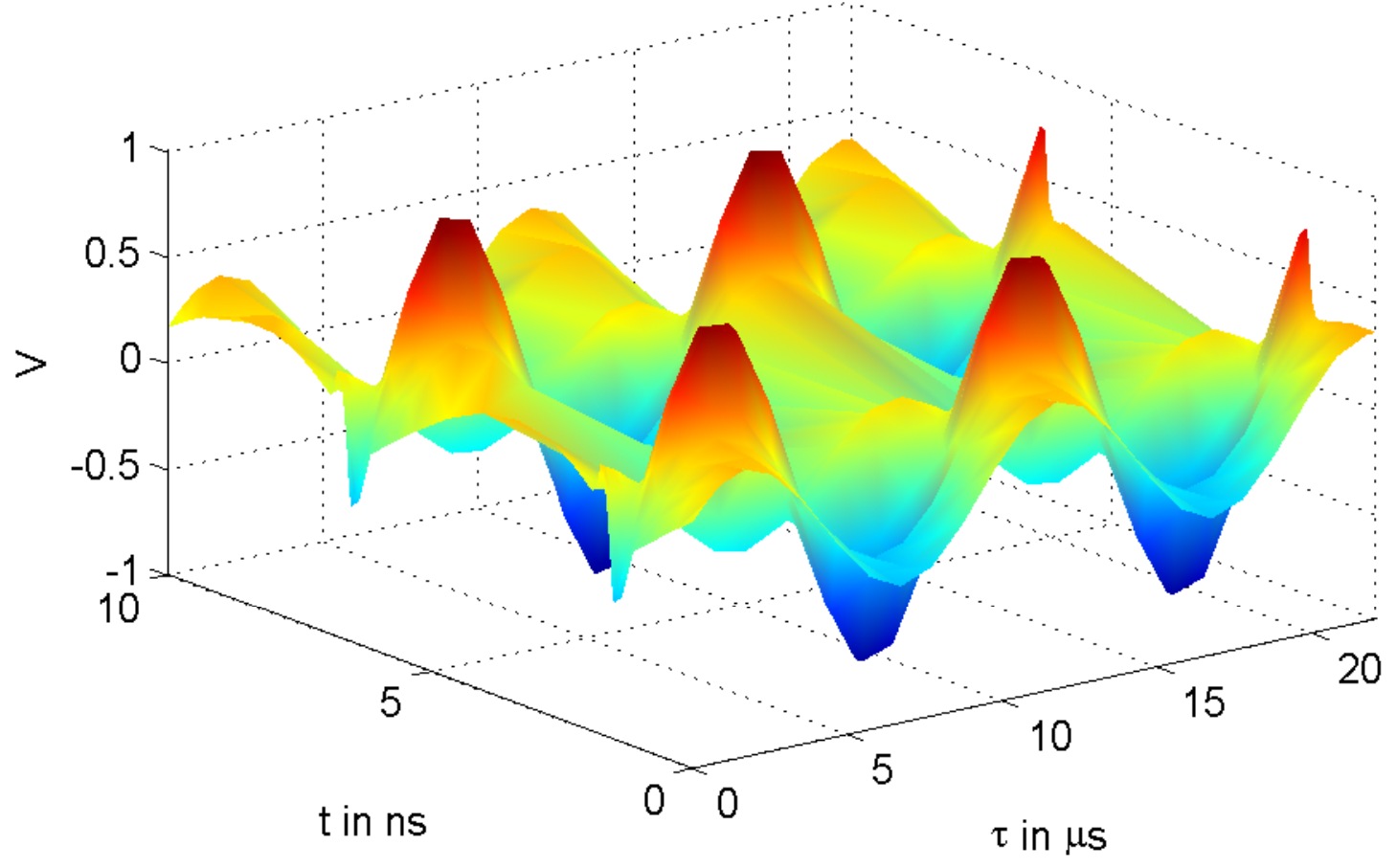}
 \fcaption{\label{mixer} Multi-rate solution for Gilbert mixer.}
 \end{center}
\end{figure}

Our first example is the simulation of a down-converting Gilbert mixer. The inputs are a sinusoidal
frequency of $f_{\mathrm{rf}}=99.9$MHz and a digital reference signal of fundamental frequency 
$f_{\mathrm{lo}}=100$MHz, which switches between $-5$V and 5V. Since the driving
frequencies are known, we have chosen $\omega(\tau)=1$ and carrier period to $P =10^{-8}\mathrm{s}$,
which corresponds to the frequency of 100MHz.
In Fig.~\ref{mixer} shows the resulting output signal.
In $\tau$-direction one can see the smooth sinusoidal envelope of frequency 
$f_{\mathrm{out}}=f_{\mathrm{lo}}-f_{\mathrm{rf}}=0.1$MHz, while the $t$ direction shows  
the influence of the carrier signal. In particular, at the jump discontinuities of the input signal
one can observe an overshoot of the output signal, due to capacitive effects.
The simulation with 28 envelope time steps for 20000 periods needed 10s, while a
transient simulation was done in 8min.

\begin{figure}
 \begin{center}
 \includegraphics[width=\columnwidth]{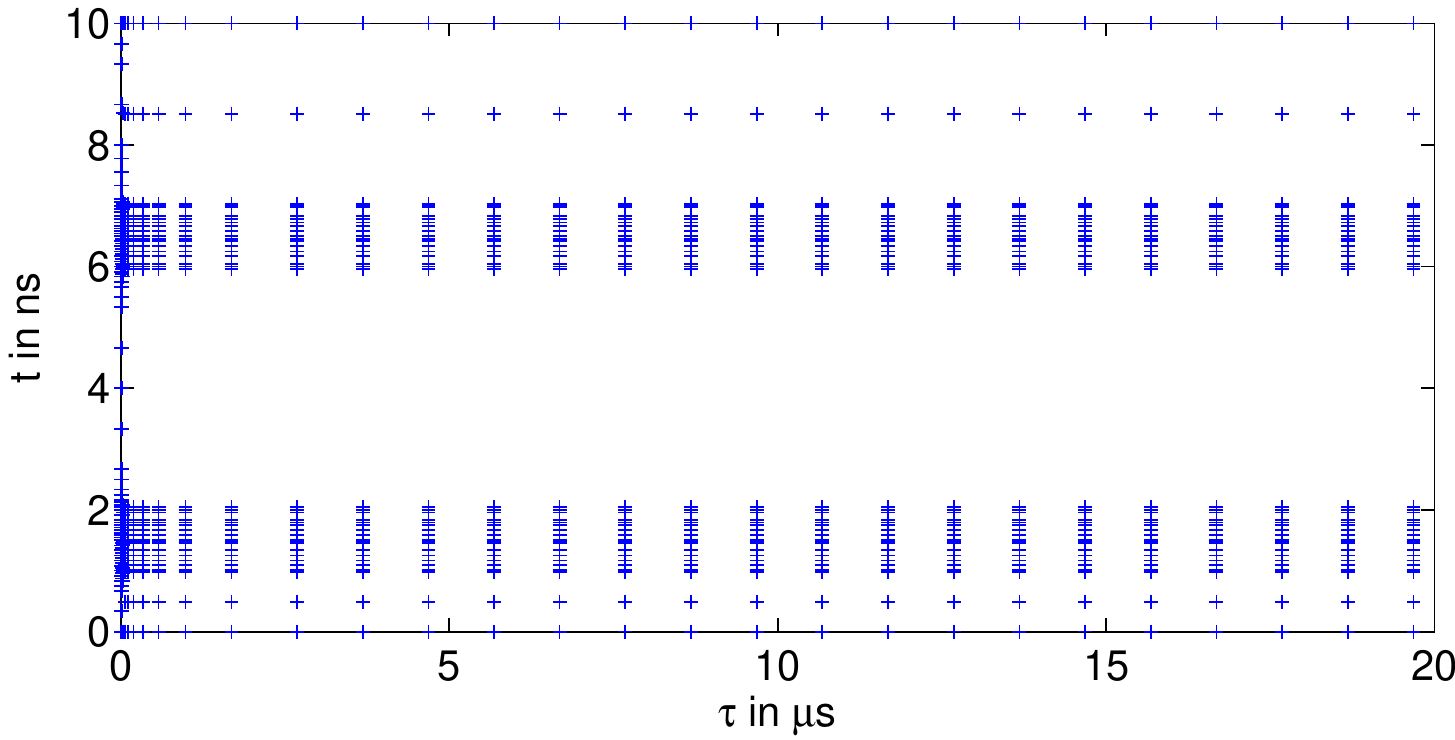}
 \fcaption{\label{mixer_grid} Adaptive grid for Gilbert mixer.}
 \end{center}
\end{figure}

In Fig.~\ref{mixer_grid} the corresponding adaptive spline grids for the envelope time steps $\tau_k$
are shown. As expected, it can be observed that the adaptive wavelet refinement leads to a fine grid
at the jump locations, while only few spline knots are needed to represent the smooth parts
in between. Furthermore we can see that roughly 30 envelope time steps are needed for an interval which 
contains 2200 oscillations of the reference signal.

Obviously the spline grid does not change much with $\tau$, since the locations of sharp transients
do not change with $\tau$. However, we will consider examples in the sequel, where an optimal grid is not known in advance,
but may even change with $\tau$.

\rengSubsection{Colpitts Oscillator}

\begin{figure}
 \begin{center}
 \includegraphics[width=\columnwidth]{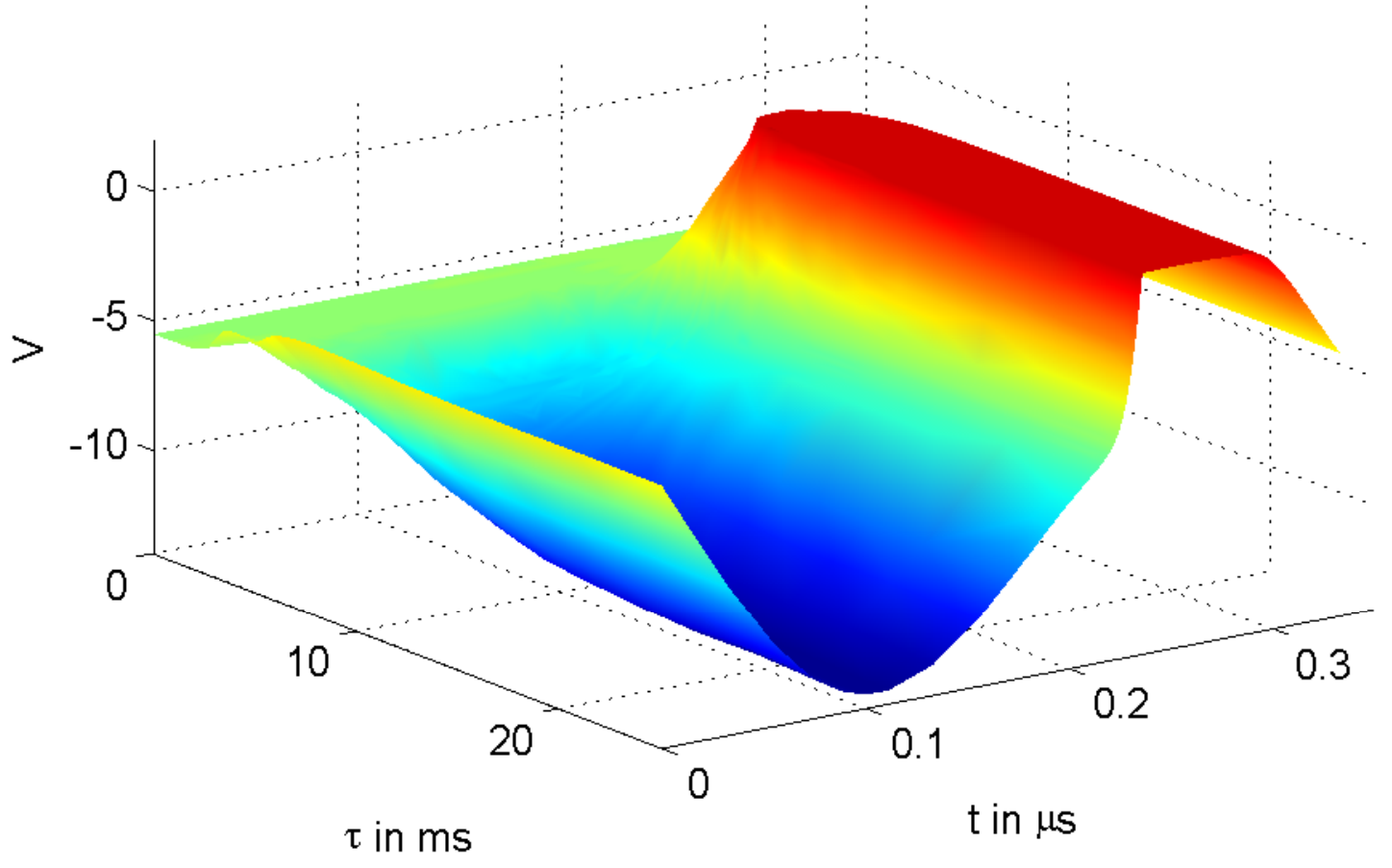}
 \fcaption{\label{oscillator} Multi-rate solution for 3MHz Colpitts oscillator.}
 \end{center}
\end{figure}

In our next example we consider the start up phase a  3MHz Colpitts-Quartz oscillator. To avoids unwanted effects from numerical
damping, we have used trigonometric splines \cite[Sect.~10.8]{Sch81} instead of the usual polynomial splines
(see \cite{BiBra12a} for details). The transistor used as feedback amplifier introduces some nonlinear effects into the output
(see Fig.~\ref{oscillator}), which results in a sharp edge near $t=0.24\mathrm{\mu s}$.
The simulation with 202 envelope time steps for 7500 periods needed 10s, while a
transient simulation was done in 165min.

\begin{figure}
 \begin{center}
 \includegraphics[width=\columnwidth]{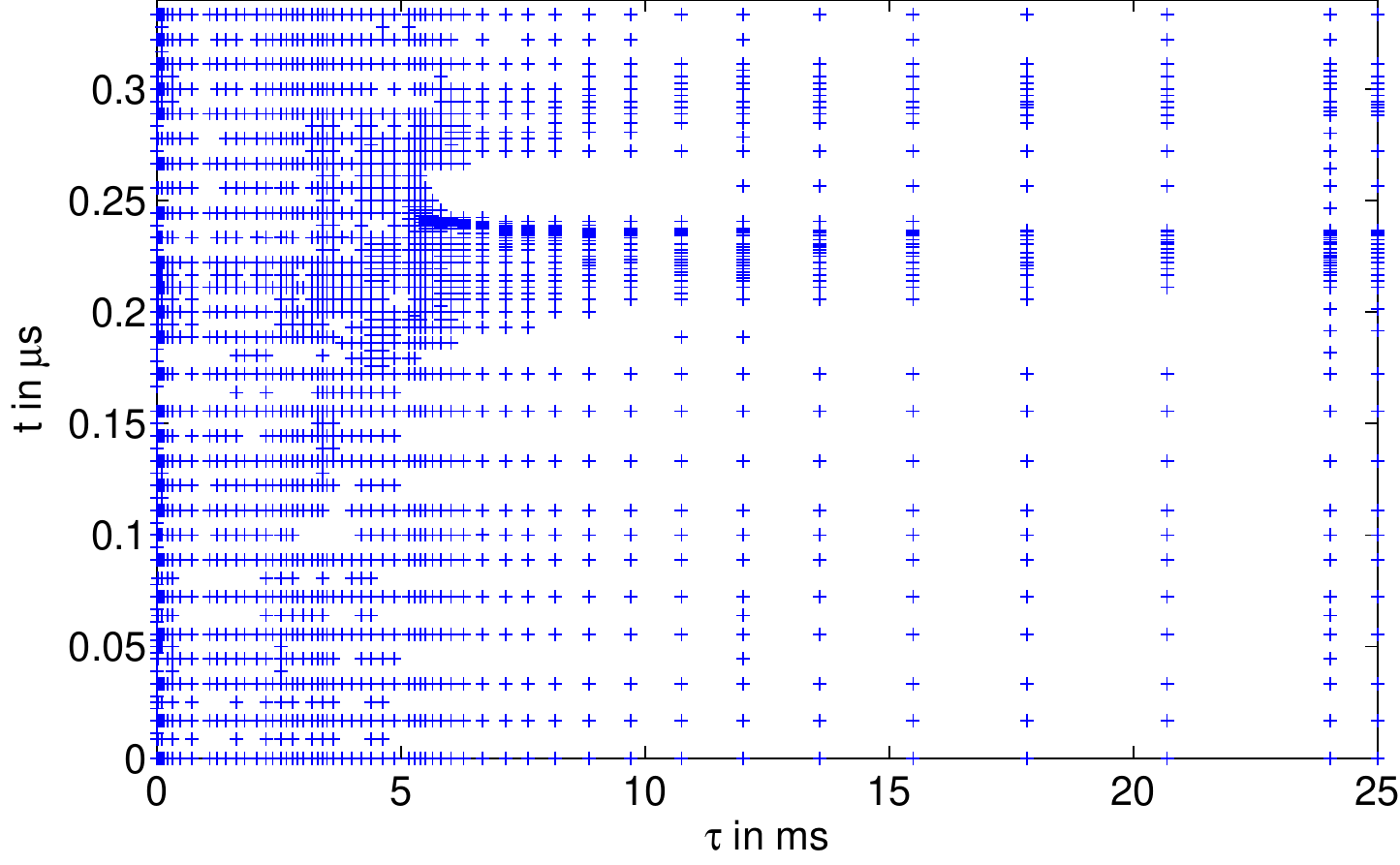}
 \fcaption{\label{oscillator_grid} Adaptive grid for Colpitts oscillator.}
 \end{center}
\end{figure}

The spline grids in Fig.~\ref{oscillator_grid} are adapted to this edge, when it becomes prominent at 
$\tau\approx 5\mathrm{ms}$, and 
the grids are even adapted to a slight change of location in the following envelope time steps, which shows the perfect 
functioning of our approach.

\rengSubsection{A Phase Locked Loop}

\begin{figure}
 \begin{center}
 \includegraphics[width=\columnwidth]{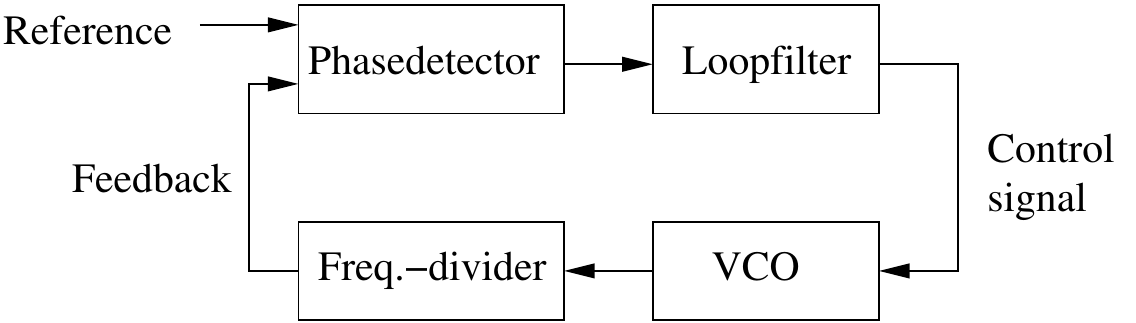}

 \fcaption{\label{pll_sch} Schematic of a PLL.}
 \end{center}
\end{figure}

Our next example is a Phase Locked Loop (PLL), 
which generates an output signal
whose phase is related to the phase of the input ``reference'' signal.
As depicted in Fig.~\ref{pll_sch}, this is achieved, by comparing the phases of the ``reference'' signal and the
output or ``feedback'' signal in a phase detector. The result is than filtered,
in order to stabilize the behavior of the PLL, and fed to a Voltage Controlled Oscillator
whose frequency depends on the (filtered) phase difference. The output is fed back
(possibly to a frequency divider) to the phase detector. If the PLL is locked, the phases
of reference and feedback signal are synchronized. Thus both signal will have (almost) 
the same local frequency.

\begin{figure}
 \begin{center}
 \includegraphics[width=\columnwidth]{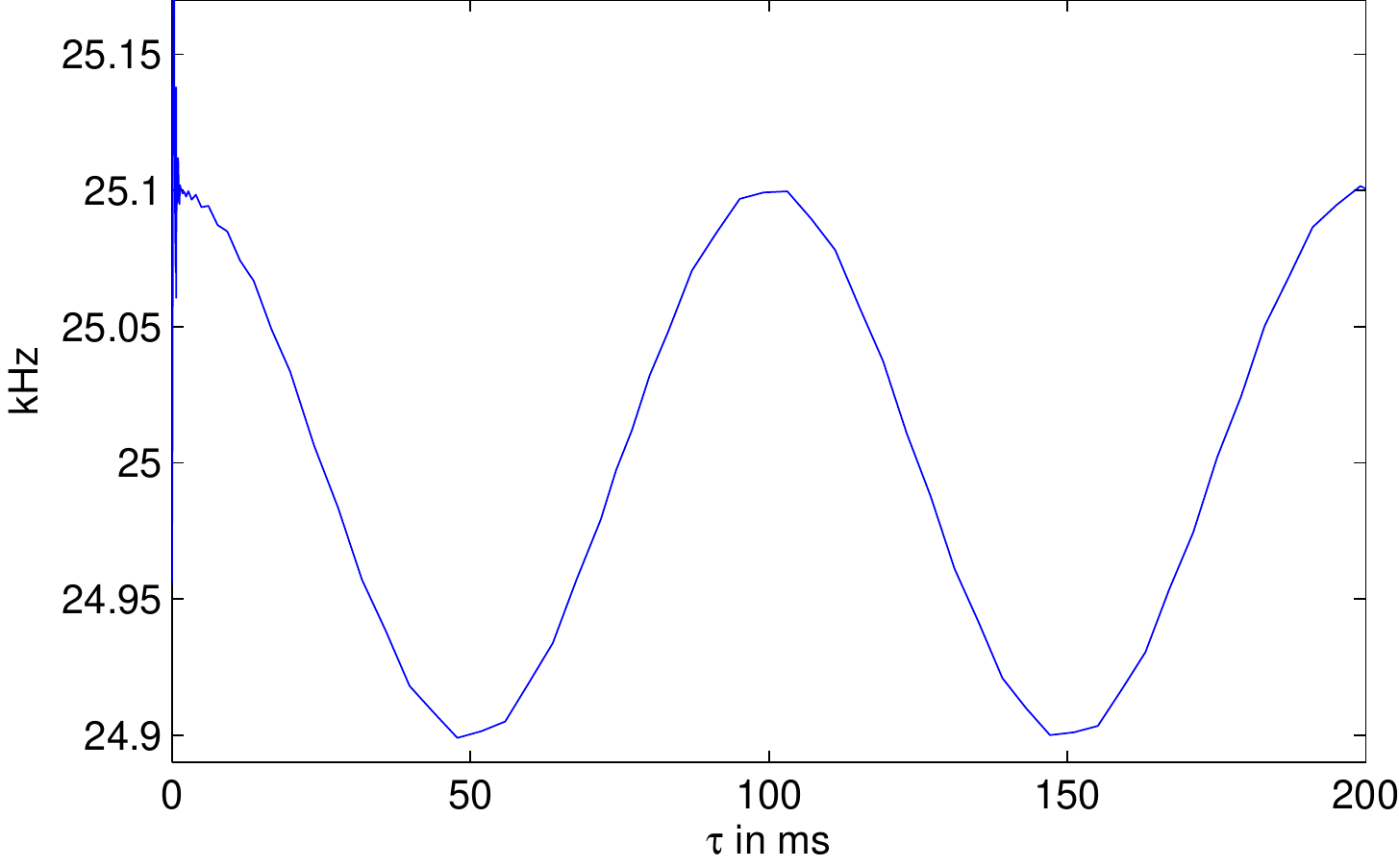}

 \fcaption{\label{pllfreq} Estimated instantaneous frequency for the PLL.}
 \end{center}
\end{figure}

We have tested a PLL (containing 205 MOSFETs) 
with a frequency modulated sinusoidal signal with center frequency 25kHz, which corresponds to 800kHz
oscillator frequency due to a frequency divider of factor 32.
The baseband signal is also sinusoidal with frequency 10Hz and frequency deviation 100Hz.
The estimated instantaneous frequency in Fig.~\ref{pllfreq} (determined based on $\omega(\tau)$)
shows after a short startup phase a good agreement with the baseband signal.

\begin{figure}
 \begin{center}
 \includegraphics[width=\columnwidth]{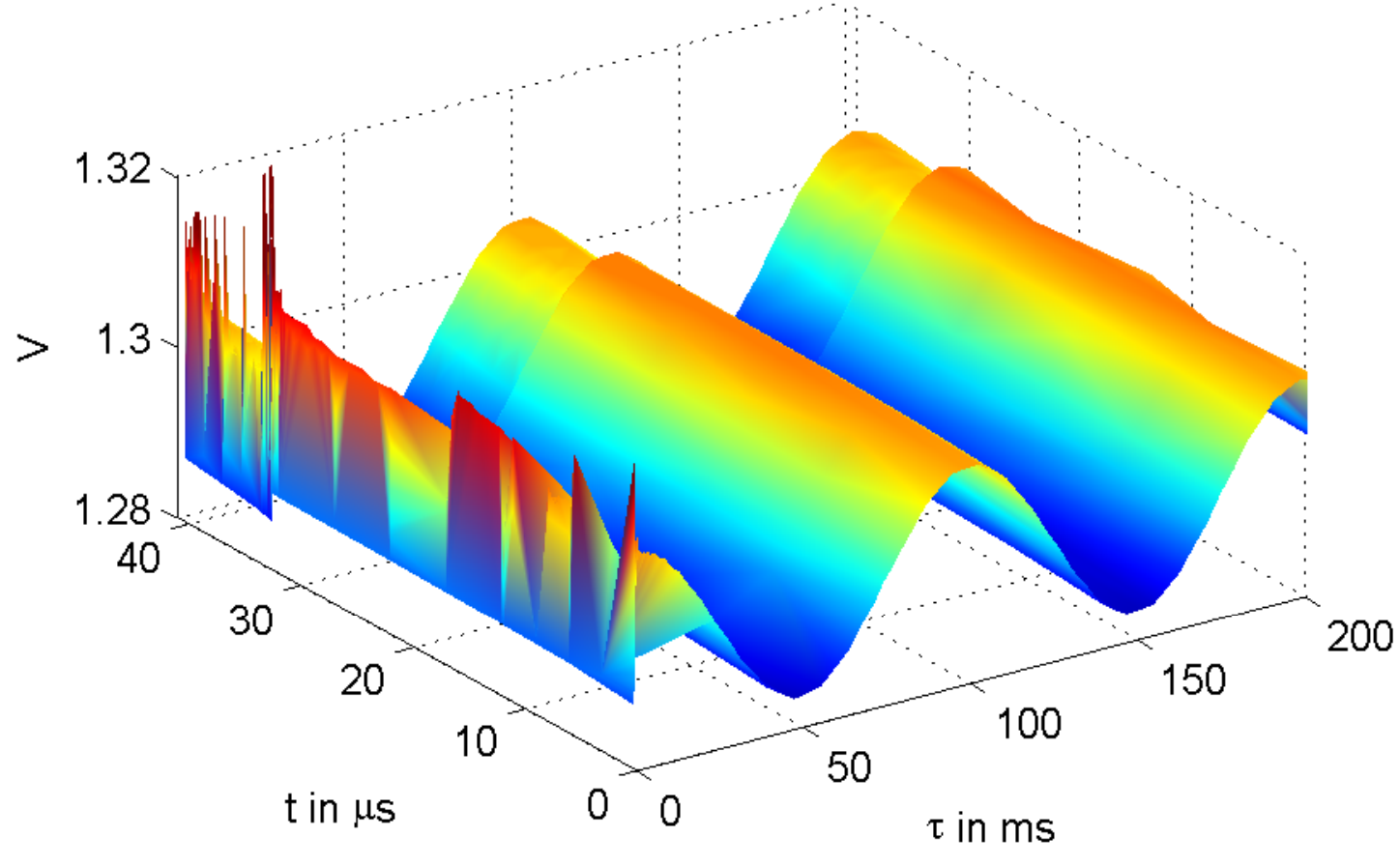}

 \fcaption{\label{pllcontrol} PLL example. Control signal of the VCO (multi-rate).}
 \end{center}
\end{figure}

\begin{figure}
 \begin{center}
 \includegraphics[width=\columnwidth]{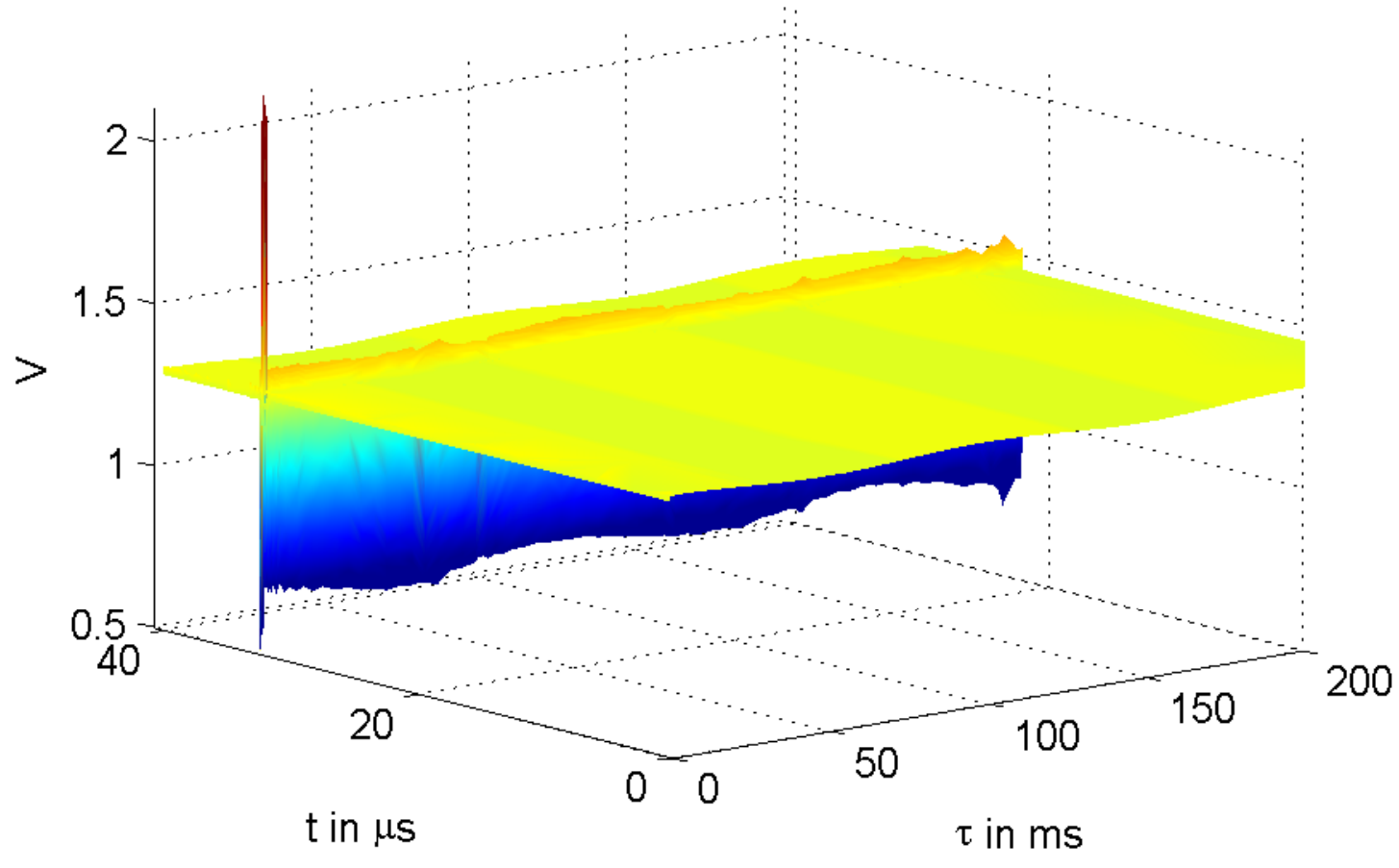}

 \fcaption{\label{pump} PLL example. Output of the phase detector (multi-rate).}
 \end{center}
\end{figure}

This is also reflected by the control signal of the VCO in Fig.~\ref{pllcontrol}, which is obtained 
by filtering the output of the phase detector (Fig.~\ref{pump}). Note, that the envelope
corresponds to the baseband signal, while the carrier signal is due to the filtering almost
constant.
This allows to do the shown envelope simulation using only 160 envelope time steps, while a 
corresponding transient analysis would contain 5000 oscillations. 
The wavelet envelope simulation of this circuit was
done in 7\,min. However, a comparable transient simulation on the same simulator needed 33\,hours. 

\begin{figure}
 \begin{center}
 \includegraphics[width=\columnwidth]{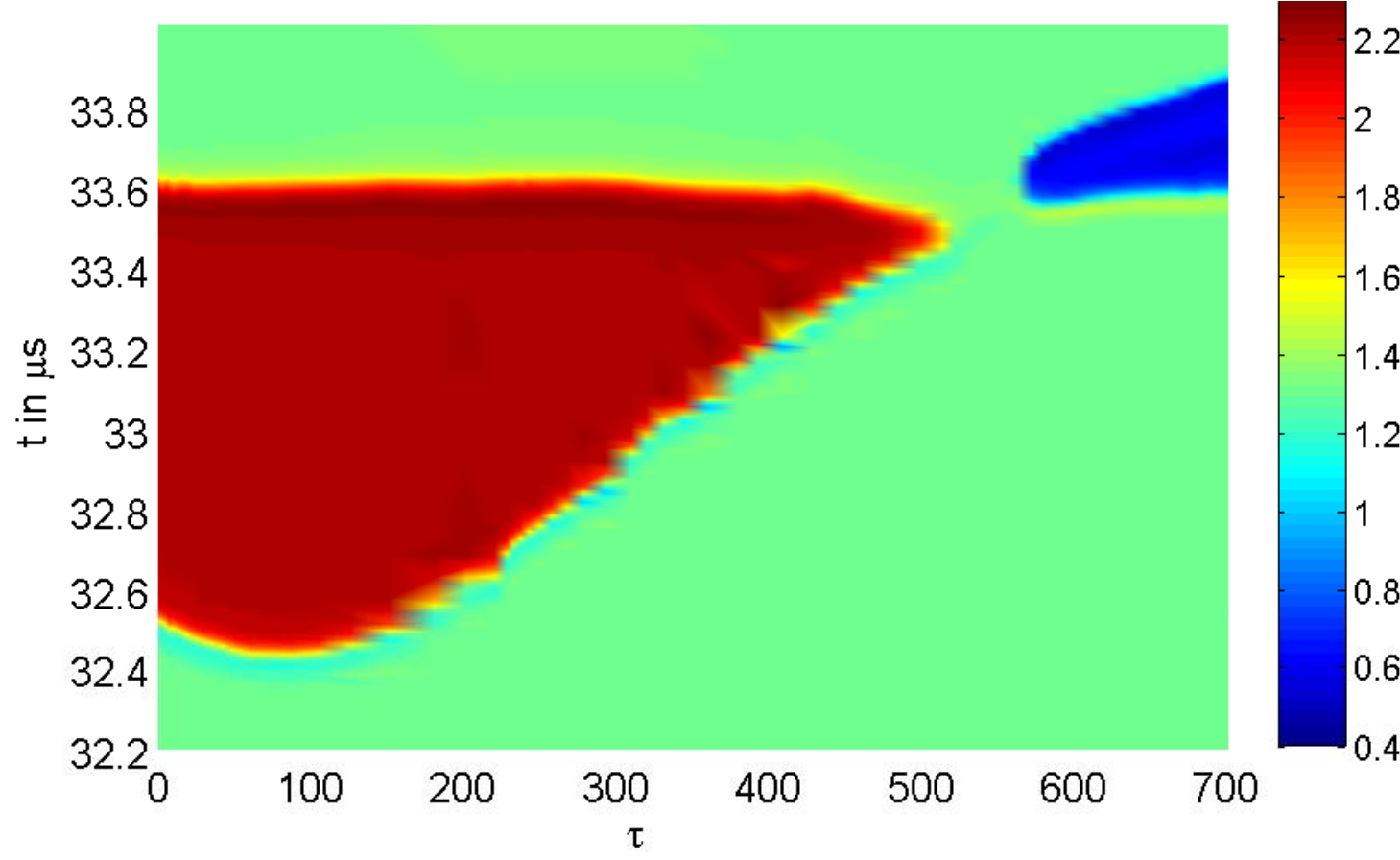}

 \fcaption{\label{plldet} PLL example. Detail of the phase detector output from Fig.~\ref{pump}.}
 \end{center}
\end{figure}

We will have a closer look to the start-up phase before the PLL is locked, which gives a good impression
how the adaptive wavelet method works. For this we have a closer look at the phase detector output for $t$ near 33ms, where 
sharp edges are present as one can see in Fig.~\ref{plldet}. The change in the signal shape is due to changes in the
phase difference (with $\tau$) during transient response.

\begin{figure}
 \begin{center}
 \includegraphics[width=\columnwidth]{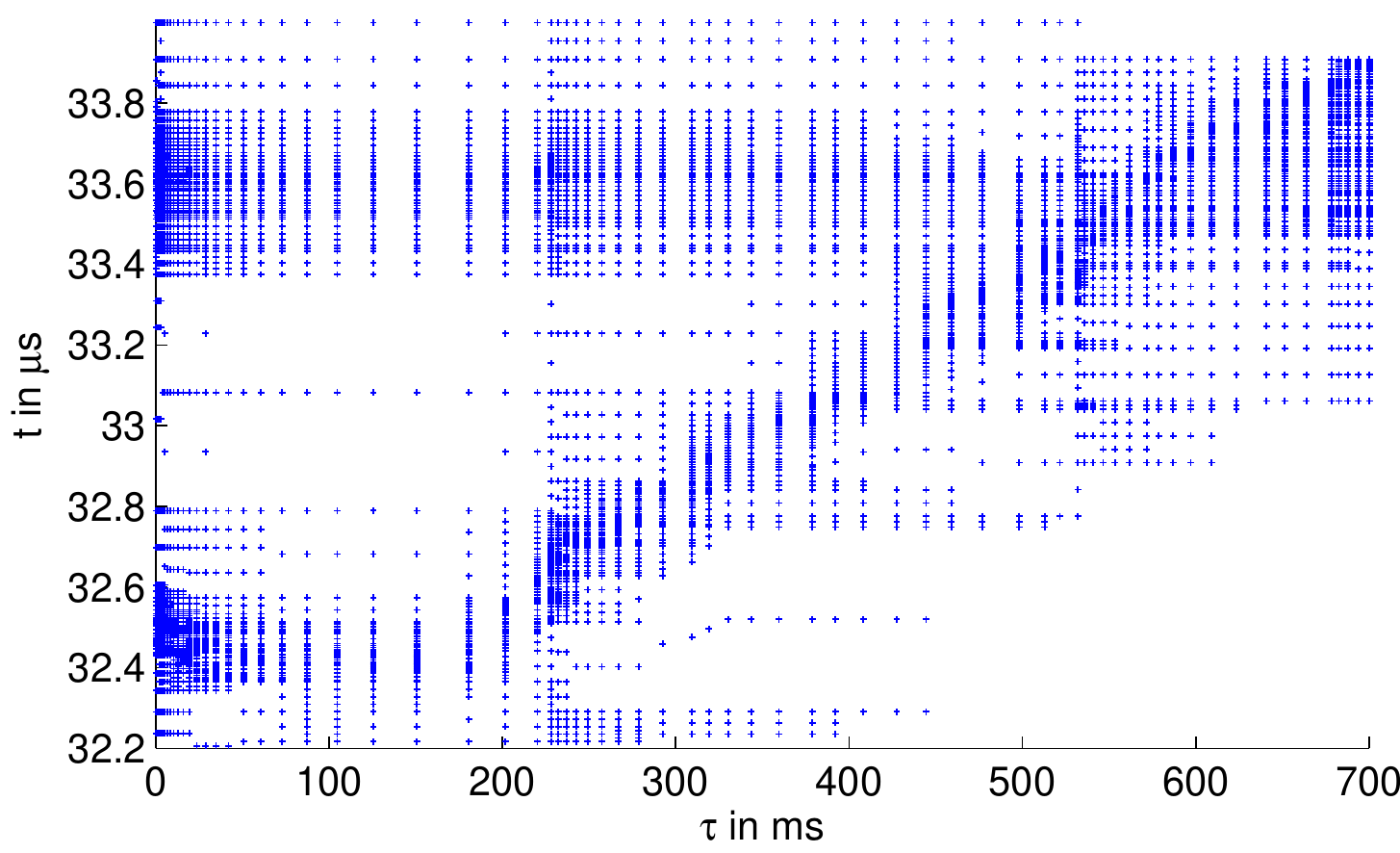}

 \fcaption{\label{pllgrid} PLL example. Spline grid corresponding to the cutout in Fig.~\ref{plldet}.}
 \end{center}
\end{figure}

In Fig.~\ref{pllgrid} one can see the corresponding spline grid. Obviously, the grid is refined around the sharp edges in the
signal, as expected

\vspace{6em}
\noindent{\Large\bfseries Acknowledgements}

This work has been 
partly supported by the ENIAC research project ARTEMOS under grant 829397 and
the FWF under grant P22549.


\vspace{0cm}
\begin{center}
\noindent{\Large\bfseries About Authors\dots}
\vspace{0cm}
\end{center}

\noindent\textbf{Kai BITTNER} received the Dipl.-Math.\  degree from the University of  
Rostock, Germany in 1996 and the PhD degree from the Technical  
University Munich, Germany, in 2000. From 2000 to 2002, he was a  
Post-Doctoral Researcher at the University of Missouri ---  
St.~Louis. In 2002, he joined the academic staff of the University  
of Ulm. From 2008 until 2010 he worked as a Researcher in the  
ICESTARS project funded by the European Union. He is currently  
employed at the University of Applied Science of Upper Austria. His  
research interests are wavelets, numerical analysis of differential  
equations, approximation theory, and circuit simulation.  

\noindent\textbf{Hans Georg BRACHTENDORF} graduated in Electrical Engineering from the RWTH Aachen, 
Germany in 1989 and obtained the Dr.-Ing.\ degree from the University of Bremen, 
Institute for Electromagnetic Theory and Microelectronics in 1994, 
also in Electrical Engineering.
From 1994-2001 he was an Assistant Professor (C1) also at the University of Bremen 
and obtained the Venia Legendi (Habilitation) 
from the same university in 2001. In 1997-1998 he was 
affiliated with the Wireless Laboratory of Bell Laboratories/Lucent Technologies in 
Murray Hill/New Jersey, where he performed research in circuit simulation and design. 
In 2001 he joined the Fraunhofer Institute for Integrated Circuits in Erlangen, Germany. 
His focus there was on system design and simulation for satellite broadcasting 
systems (WorldSpace, XM radio) and transceiver designs.
Since 2005 Dr. Brachtendorf is full professor at the Fachhochschule Ober\"osterreich 
for System Design and Simulation, Communications and Signal Processing. 
He is author and co-author of 1 book and numerous technical papers 
dealing mainly with circuit simulation and device modeling. 
He holds four patents in various fields of circuit analysis and design, 
including a patent for a novel image reject filter, a subsampling receiver 
architecture and on multirate simulation techniques.
His research interests encompass circuit design, modeling and simulation 
as well as signal processing and digital communication. 

\end{multicols}
\end{document}